\documentclass{amsart}

\usepackage{amssymb}
\usepackage{amsmath}
\usepackage{amsthm}
\usepackage{amsfonts}

\usepackage{latexsym}

\usepackage[dvips]{epsfig}
\usepackage{amsbsy}
\usepackage{amsgen}
\usepackage{amscd}
\usepackage{amsopn}
\usepackage{amstext}
\usepackage{amsxtra}
\usepackage{xypic}

\newtheorem{theorem}{Theorem}[section]

\newtheorem{examples}[theorem]{Examples}
\newtheorem{remark}[theorem]{Remark}
\newtheorem{lemma}[theorem]{Lemma}
\newtheorem{proposition}[theorem]{Proposition}

\newtheorem{notation}[theorem]{Notation}

\newtheorem{calculation}[theorem]{Calculation}
\newtheorem{fact}[theorem]{Fact}

\begin{document}

\title[Formal aspects of Gray's tensor products of 2-categories]
{Formal aspects of Gray's tensor product of 2-categories}
\author[{A. E.} {Stanculescu}]{{Alexandru E.} {Stanculescu}}
\address{\newline Department of Mathematics and Statistics,
\newline Masaryk University,  Kotl\'{a}{\v{r}}sk{\'{a}} 2,\newline
611 37 Brno, Czech Republic}
\email{stanculescu@math.muni.cz}
\thanks{Research supported by the Ministry of Education of the 
Czech Republic under grant LC505}

\begin{abstract}
The category of small $2$-categories has two monoidal structures 
due to John Gray: one biclosed and one closed. We propose a 
formalisation of the construction of the right internal 
and internal homs of these monoidal structures.
\end{abstract}
\maketitle
\section{Introduction}

Let {\bf Cat} be the category of small categories.  In \cite{Gr1}, Gray
introduced two monoidal structures on the category $2$-{\bf Cat} of 
small categories enriched over {\bf Cat}: one biclosed and one closed.  
The coherence axioms for these monoidal structures were proved 
in \cite{Gr2}; see also \cite{Bou}. We recall that the right internal hom
of the biclosed monoidal structure on $2$-{\bf Cat} consists of
$2$-functors, quasi-natural (also called lax natural) transformations
and modifications.

This paper is an attempt to understand these 
monoidal structures and some facts surrounding them.
In this respect we propose a formalisation of the 
construction of the (right) internal 
hom of $2$-{\bf Cat}, as outlined below.

Let $\mathcal{V}$ be a closed category and let 
$\mathcal{V}$\text{-}{\bf Cat} be the category of 
small $\mathcal{V}$-categories. Given a comonoid
$C^{\bullet}$ in the category of cosimplicial objects in $\mathcal{V}$
equipped with the Day convolution product, and two small 
$\mathcal{V}$-categories $\mathcal{A}$ and $\mathcal{B}$, 
we construct a small $\mathcal{V}$-category 
$Coh^{C^{\bullet}}(\mathcal{A},\mathcal{B})$ whose 
objects are the $\mathcal{V}$-functors from $\mathcal{A}$
to $\mathcal{B}$ and whose homs are the objects of 
$\mathcal{V}$-coherent transformations with respect to $C^{\bullet}$.
The notation $Coh$ is borrowed from the work of Cordier 
and Porter \cite{CP}, who made a similar construction when
$\mathcal{V}$ is the category of simplicial sets, except that in 
their case $Coh(\mathcal{A},\mathcal{B})$ is not a category
enriched over simplicial sets. When $C^{\bullet}$ is the constant 
cosimplicial object with value the unit object
of $\mathcal{V}$, we recover the standard internal 
hom of $\mathcal{V}$\text{-}{\bf Cat} consisting of the 
$\mathcal{V}$-natural transformations. When 
$\mathcal{V}$={\bf Cat} and $C^{\bullet}$ is what we call 
the standard cocategory (cogroupoid) interval in {\bf Cat}, 
our construction recovers the right internal (internal) hom of 
$2$-{\bf Cat}. We actually construct a $\mathcal{V}$-category 
$Coh^{C}(\mathcal{A},\mathcal{B})$ for every 
comonoid $C$ in the category of coaugmented 
cosimplicial objects in $\mathcal{V}$.
The construction is natural in all three variables, and we
present it decomposed in as many steps as we could. For 
instance, the endofunctor $Coh^{C}(\mathcal{A},-)$ is a 
composite of four natural functors, each one having left adjoints.

As an outcome of this formalisation we obtain a
formula for Gray's tensor product(s) of $2$-categories. 
Other outcomes will be detailed in \cite{St}.

The paper is organised as follows. Sections 2 and 3 are preparatory,
and consist of recollections of facts regarding the Day convolution
products on certain functor categories, actions of monoidal categories
and Tensor-Hom situations. In section 4 we exhibit a chain 
of monoidal functors which will be part of the construction 
of $Coh^{C}(\mathcal{A},-)$. The most important one is a 
familiar cosimplicial cobar construction. In section 5 we construct 
$Coh^{C}(\mathcal{A},\mathcal{B})$ and show that, as a functor 
of three variables, it is part of a Tensor-Hom situation. We end by 
showing that if we take $\mathcal{V}$={\bf Cat}, equipped with the 
cartesian closed structure, then the standard cocategory interval in 
{\bf Cat}, to be denoted by $\mathbb{I}^{\bullet}$, is a comonoid in 
the category of cosimplicial objects in {\bf Cat}, and that 
$Coh^{\mathbb{I}^{\bullet}}(\mathcal{A},\mathcal{B})$ is 
precisely the right internal hom of $2$-{\bf Cat}.

\section{Background, part one}
We denote by $\Delta$ the category of finite 
non-empty ordinals and order preserving maps. The ordinal
$n+1=\{0,...,n\}$ will be denoted by $[n]$. We let $\Delta_{+}$ 
be the category of all finite ordinals and order preserving maps. 
The ordinal $n$ will be denoted by $n$. $\Delta_{+}$ has a monoidal 
product given by the ordinal addition, with unit the ordinal $0$.
We denote by $i$ the inclusion $\Delta \subset \Delta_{+}$. We 
let $\Delta(n)$ be the $n$-th truncation of $\Delta$ and 
$\Delta(n)_{mon}$ be $\Delta(n)$ without the 
codegeneracies $s^{i}:[n]\rightarrow [n-1]$, $n\geq 1$.
If $\mathcal{V}$ is a monoidal category, we denote by 
$Comon(\mathcal{V})$ the category of comonoids in $\mathcal{V}$.

Throughout this section $(\mathcal{V},\otimes, I)$ is a cocomplete 
closed category.
\subsection{Coreflexive graphs}
The category $\mathcal{V}^{\Delta(1)}$ admits the non-symmetric 
Day convolution product $\bigstar$. We recall its construction in detail.
For $X^{\bullet},Y^{\bullet}\in \mathcal{V}^{\Delta(1)}$, one has
$(X^{\bullet}\bigstar Y^{\bullet})^{0}=X^{0}\otimes Y^{0}$
and $(X^{\bullet}\bigstar Y^{\bullet})^{1}$ is the pushout of the diagram
\[
   \xymatrix{
X^{0}\otimes Y^{1} \ar[rr]^{i_{X^0,Y^1}} & & 
(X^{\bullet}\bigstar Y^{\bullet})^{1}\\
X^{0}\otimes Y^{0} \ar[rr]^{d^{1}\otimes 1_{Y^{0}}} 
\ar[u]^{1_{X^{0}}\otimes d^{0}} & & X^{1}\otimes Y^{0} \ar[u]_{i_{X^1,Y^0}}\\
}
  \]
The unit of $\bigstar$ is $cstI$, the constant 
1-truncated cosimplicial object with value $I$.
The cofaces are $D^{0}=i_{X^1,Y^0}(d^{0}\otimes 1_{Y^{0}})$ 
and $D^{1}=i_{X^0,Y^1}(1_{X^{0}}\otimes d^{1})$. 
The codegeneracy is obtained using the universal property of the pushout. 
The associativity isomorphism can be seen from the diagram
\[
   \xymatrix{
X^{0}\otimes Y^{0}\otimes Z^{0} \ar @{..>} [rr]^{d^{1}\otimes 1_{Z^{0}}} 
\ar[dr]_{1_{X^{0}}\otimes d^{0}}
\ar @<-4pt> @{..>} [dd]_{1_{X^{0}}\otimes d^{0}\otimes 1_{Z^{0}}} 
\ar @<4pt> [dd]^{1_{X^{0}}\otimes d^{1}\otimes 1_{Z^{0}}}
& & X^{1}\otimes Y^{0}\otimes Z^{0} \ar @{..>} [dd]\\
&  X^{0}\otimes Y^{0}\otimes Z^{1}\ar[dd]\\
X^{0}\otimes Y^{1}\otimes Z^{0} \ar @{..>} [rr]|{i_{X^0,Y^1}\otimes 1_{Z^{0}}}
\ar[dr]_{1_{X^{0}}\otimes i_{Y^1,Z^0}} & & 
(X^{\bullet}\bigstar Y^{\bullet})^{1}\otimes Z^{0} \ar @{-->} [dr]\\
&  X^{0}\otimes (Y^{\bullet}\bigstar Z^{\bullet})^{1} \ar @{-->} [rr] & &
(X^{\bullet}\bigstar Y^{\bullet} \bigstar Z^{\bullet})^{1}\\
 }
   \]
where all the faces are pushouts. The object 
$(X^{\bullet}\bigstar (Y^{\bullet} \bigstar Z^{\bullet}))^{1}$
is obtained from the back face and the map $1_{X^{0}}\otimes i_{Y^1,Z^0}$; 
the object $((X^{\bullet}\bigstar Y^{\bullet}) \bigstar Z^{\bullet})^{1}$ is 
obtained from the left face and the map $i_{X^0,Y^1}\otimes 1_{Z^{0}}$.

The monoidal product $\bigstar$ restricts 
to a monoidal product $\bigstar$ on $\mathcal{V}^{\Delta(1)_{ mon}}$.

\subsection{(Coaugmented) cosimplicial objects} 
The category $\mathcal{V}^{\Delta}$ has the non-symmetric 
Day convolution product $\bigstar$, see \cite{Ba}, \cite{CP}, \cite{MS}. 
We detail two of its presentations. If 
$X^{\bullet},Y^{\bullet}\in \mathcal{V}^{\Delta}$,
one has $(X^{\bullet}\bigstar Y^{\bullet})^{0}=X^{0}\otimes Y^{0}$.
For $n\geq 1$, $(X^{\bullet}\bigstar Y^{\bullet})^{n}$ is the coequaliser
\[
\xymatrix{
\underset{p+q=n-1}\coprod X^{p}\otimes Y^{q} 
\ar @<0pt> [rr]_{v} \ar @<6pt> [rr]^{u} & &
\underset{r+s=n}\coprod X^{r}\otimes Y^{s}\\
X^{p}\otimes Y^{q} \ar[rr]^{1_{X^{p}}
\otimes d^{0}} \ar[u]_{\mathrm{inj}^{n-1}_{p,q}}
\ar[dr]_{d^{p+1}\otimes 1_{Y^{q}}} & &
X^{p}\otimes Y^{q+1} \ar[u]_{\mathrm{inj}^{n}_{p,q+1}}\\
 & X^{p+1}\otimes Y^{q} \ar[uur]_{\mathrm{inj}^{n}_{p+1,q}}\\
  }
  \]
where $u\mathrm{inj}^{n-1}_{p,q}=
\mathrm{inj}^{n}_{p+1,q}(d^{p+1}\otimes 1_{Y^{q}})$ 
and $v\mathrm{inj}^{n-1}_{p,q}=
\mathrm{inj}^{n}_{p,q+1}(1_{X^{p}}\otimes d^{0})$. 
For $0\leq k\leq n+1$, the coface map $D^{k}$
is determined by the diagram
\[
\xymatrix{
X^{r}\otimes Y^{s-1} \ar[rr]^{1_{X^{r}}\otimes d^{0}} 
\ar[dr]^{\mathrm{inj}^{n-1}_{r,s-1}} \ar[dd]^{1_{X^{r}}\otimes d^{k-r-1}}
& & X^{r}\otimes Y^{s} \ar[dd]^{1_{X^{r}}\otimes d^{k-r}} 
\ar[dr]^{\mathrm{inj}^{n}_{r,s}}\\
&  \underset{p+q=n-1}\coprod X^{p}\otimes Y^{q} \ar[rr]^{v} 
& & \underset{r+s=n}\coprod X^{r}\otimes Y^{s}\\
X^{r}\otimes Y^{s} \ar[rr]|{1_{X^{r}}\otimes d^{0}} 
\ar[dr]^{\mathrm{inj}^{n}_{r,s}} & &
X^{r}\otimes Y^{s+1} \ar[dr]^{\mathrm{inj}^{n+1}_{r,s+1}}\\
& \underset{p+q=n}\coprod X^{p}\otimes Y^{q} \ar[rr]|{v} 
& & \underset{r+s=n+1}\coprod X^{r}\otimes Y^{s}\\
 }
   \]
if $r< k$, and by the diagram
\[
\xymatrix{
X^{r-1}\otimes Y^{s} \ar[rr]^{d^{r}\otimes 1_{Y^{s}}} 
\ar[dr]^{\mathrm{inj}^{n-1}_{r-1,s}} \ar[dd]^{d^{k}\otimes 1_{Y^{s}}}
& & X^{r}\otimes Y^{s} \ar[dd]^{d^{k}\otimes 1_{Y^{s}}} 
\ar[dr]^{\mathrm{inj}^{n}_{r,s}}\\
&  \underset{p+q=n-1}\coprod X^{p}\otimes Y^{q} \ar[rr]^{u} & & 
\underset{r+s=n}\coprod X^{r}\otimes Y^{s}\\
X^{r}\otimes Y^{s} \ar[rr]|{d^{r+1}\otimes 1_{Y^{s}}} 
\ar[dr]^{\mathrm{inj}^{n}_{r,s}} & &
X^{r+1}\otimes Y^{s} \ar[dr]^{\mathrm{inj}^{n+1}_{r+1,s}}\\
& \underset{p+q=n}\coprod X^{p}\otimes Y^{q} \ar[rr]|{u} & &
\underset{r+s=n+1}\coprod X^{r}\otimes Y^{s}\\
 }
   \]
if $r\geq k$. The codegeneracies are defined similarly,
see \cite{Ba},\cite{MS}. The unit of $\bigstar$ is $cstI$, the 
constant cosimplicial object with value $I$. 
$(X^{\bullet}\bigstar Y^{\bullet})^{n}$ $(n\geq 1)$ 
can also be calculated \cite{Ba} as the colimit of the diagram
\[
   \xymatrix{
X^{0}\otimes Y^{n}\\
X^{0}\otimes Y^{n-1} \ar[u]^{1_{X^{0}}\otimes d^{0}} 
\ar[r]^{d^{1}\otimes 1_{Y^{n-1}}} & X^{1}\otimes Y^{n-1}\\
& ...\ar[u]\\
& \ar[r] & X^{p}\otimes Y^{q+1}\\
& & X^{p}\otimes Y^{q} \ar[u]^{1_{X^{p}}\otimes d^{0}} 
\ar[r]^{d^{p+1}\otimes 1_{Y^{q}}} & X^{p+1}\otimes Y^{q}\\
& & & \ar[u]\\
& & & ...\\
& & & X^{n-1}\otimes Y^{0} \ar[u]^{1_{X^{n-1}}\otimes d^{0}} 
\ar[r]^{d^{n}\otimes 1_{Y^{0}}} & X^{n}\otimes Y^{0}\\
}
   \]
where $p+q=n-1$. It is this presentation that we shall use the most.

Any monoidal (resp. opmonoidal and cocontinuous) functor
$\mathcal{V}_{1}\rightarrow \mathcal{V}_{2}$ between cocomplete 
closed categories induces a monoidal (resp. opmonoidal) functor
$\mathcal{V}_{1}^{\Delta}\rightarrow \mathcal{V}_{2}^{\Delta}$.
There are various adjoint pairs between $\mathcal{V}$ and 
$\mathcal{V}^{\Delta}$, which we summarise as 
$$sk \dashv ev_{0}\dashv cst\dashv H^{0}$$
Here $cst$ denotes the constant cosimplicial object functor, 
$ev_{0}$ is the evaluation at $[0]$ and 
$sk(A)^{n}=\underset{\Delta([0],[n])}\sqcup A$. 
The functors $cst$ and  $ev_{0}$ are strong monoidal,
and $sk$ is opmonoidal for formal reasons.
We have induced adjoint pairs
$$sk:Comon(\mathcal{V})\rightleftarrows 
Comon(\mathcal{V}^{\Delta}):ev_{0}
\ \text{and} \ ev_{0}:Comon(\mathcal{V}^{\Delta})
\rightleftarrows Comon(\mathcal{V}):cst$$
The functor $ev_{0}:Comon(\mathcal{V}^{\Delta})\rightarrow 
Comon(\mathcal{V})$ is a (Grothendieck) bifibration, provided that
$\mathcal{V}$ is sufficiently complete. The same adjunctions and the
same facts concerning the two categories of comonoids hold
if one replaces $\Delta$ by $\Delta(1)$.
\\

The functor category $\mathcal{V}^{\Delta_{+}}$ has the Day 
convolution product $\bigstar$. Its unit is $F\Delta_{+}(0,-)$,
where $F:Set\rightarrow \mathcal{V}$ is $F(S)=\underset{S}\sqcup I$.
The functor $i^{\ast}:\mathcal{V}^{\Delta_{+}}
\rightarrow \mathcal{V}^{\Delta}$ is strong monoidal. 
(To see this it suffices \cite[Theorem 5.1]{IK} 
to show that $\Delta([0],-)$ is a monoid
in $(Set^{\Delta},\bigstar)^{op}$, that is, a comonoid in 
$(Set^{\Delta},\bigstar)$. But $\Delta([0],-)=sk(1)$.)
Therefore $i_{!}$ is opmonoidal for formal reasons.
We have an induced adjoint pair $$i_{!}:Comon(\mathcal{V}^{\Delta})
\rightleftarrows Comon(\mathcal{V}^{\Delta_{+}}):i^{\ast}$$

\section{Background, part two}
Let $(\mathcal{V},\otimes,I)$ be a monoidal category.
We recall that an {\bf action} of $\mathcal{V}$ on a category $\mathcal{E}$
is the data consisting of a functor $\ast:\mathcal{V}\times \mathcal{E}
\rightarrow \mathcal{E}$ and natural isomorphisms $\alpha$ and
$\lambda$ with components $\alpha_{A,B,X}:(A\otimes B)\ast X
\rightarrow A\ast(B\ast X)$ and $\lambda_{X}:I\ast X\rightarrow X$,
subject to certain coherence conditions (see, for example, \cite{JK}).
\\

Let $\mathcal{E}_{i}$ $(1\leq i\leq 3)$ be a category.
We recall from \cite{Gr3} that a {\bf TH-situation} 
consists of two functors
$$T:\mathcal{E}_{1}\times \mathcal{E}_{2}
\rightarrow \mathcal{E}_{3}$$
$$H:\mathcal{E}_{2}^{op}\times \mathcal{E}_{3}
\rightarrow \mathcal{E}_{1}$$
and natural isomorphisms
$$\mathcal{E}_{3}(T(X_{1},X_{2}),X_{3})\cong
\mathcal{E}_{1}(X_{1},H(X_{2},X_{3}))$$

If $$\ast:\mathcal{V}\times \mathcal{W}\rightarrow \mathcal{W}$$
$$H:\mathcal{W}^{op}\times \mathcal{W}\rightarrow \mathcal{V}$$
is a TH-situation with $\ast$ an action of $\mathcal{V}$ on 
a category $\mathcal{W}$, then $\mathcal{W}$ becomes a 
tensored $\mathcal{V}$-category.
If, in addition, $\mathcal{W}$ is a monoidal category and $\ast$ is a 
strong monoidal fuctor, then $\mathcal{W}$ becomes a monoidal
$\mathcal{V}$-category. In this case, if $C$ is a comonoid in 
$\mathcal{W}$, then $\mathcal{W}(C,-):\mathcal{W}\rightarrow 
\mathcal{V}$ is a monoidal $\mathcal{V}$-functor. We shall
use these observations for $\mathcal{W}\in \{\mathcal{V}^{\Delta(1)},
\mathcal{V}^{\Delta},\mathcal{V}^{\Delta_{+}}\}$, with the
obvious action of $\mathcal{V}$ and with the 
monoidal products described in section 2.
\\

Let $\mathcal{E}_{i}$ $(1\leq i\leq 3)$ be category and let 
$$T:\mathcal{E}_{1}\times \mathcal{E}_{2}\rightarrow \mathcal{E}_{3}$$
$$H:\mathcal{E}_{2}^{op}\times \mathcal{E}_{3}\rightarrow \mathcal{E}_{1}$$
be a TH-situation. Every small category $\mathbb{C}$ 
induces an obvious TH-situation $$T^{\mathbb{C}}:
\mathcal{E}_{1}^{\mathbb{C}}\times 
\mathcal{E}_{2}^{\mathbb{C}^{op}}\rightarrow \mathcal{E}_{3}$$
$$H^{\mathbb{C}}:(\mathcal{E}_{2}^{\mathbb{C}^{op}})^{op}
\times \mathcal{E}_{3}\rightarrow \mathcal{E}_{1}^{\mathbb{C}}$$
Suppose that $\mathcal{E}_{i}$ $(1\leq i\leq 3)$ is a monoidal category, 
$T$ is a monoidal functor and all the functor categories in the preceding
TH-situation admit a Day convolution product. Then $T^{\mathbb{C}}$
is a monoidal functor. In particular, if $A$ is a monoid in 
$\mathcal{E}_{2}^{\mathbb{C}^{op}}$, the functor 
$T^{\mathbb{C}}(-,A): \mathcal{E}_{1}^{\mathbb{C}}\rightarrow 
\mathcal{E}_{3}$ is monoidal.

\section{Some monoidal functors}

Let $(\mathcal{V},\otimes,I)$ be an arbitrary monoidal category. We denote by
$\mathcal{V}$-{\bf Cat} (resp. $\mathcal{V}$-{\bf CAT})
the category of small (resp. large) $\mathcal{V}$-categories and by $Ob$
the functor sending an $\mathcal{V}$-category to its set of objects.
For a set $S$, we denote by $\mathcal{V}$-{\bf Cat}$(S)$ the category of 
small $\mathcal{V}$-categories with fixed set of objects $S$.
When $\mathcal{V}$ is symmetric monoidal, $\mathcal{V}$-{\bf Cat} is 
a symmetric monoidal category with monoidal product $\otimes$ and 
unit $\mathcal{I}$, where $\mathcal{I}$ has a single object $\ast$
and $\mathcal{I}(\ast,\ast)=I$. A monoidal functor $F:\mathcal{V}\rightarrow \mathcal{W}$ 
between monoidal categories induces a functor $F:\mathcal{V}\text{-}{\bf Cat}\rightarrow
\mathcal{W}\text{-}{\bf Cat}$. Let $\mathcal{A}$ be a $\mathcal{V}$-category. 
We denote by $\mathcal{A}^{op}$ the opposite of $\mathcal{A}$. Suppose that 
$\mathcal{V}$ is a closed category. We write $Y^{X}$
for the internal hom of two objects $X$, $Y$ of $\mathcal{V}$.
Given $\mathcal{V}$-categories $\mathcal{A}$ and $\mathcal{B}$,
we denote by $\mathcal{V}\text{-}{\bf Mod}(\mathcal{A},\mathcal{B})$
the $\mathcal{V}$-category of $\mathcal{V}$-functors 
$\mathcal{B}^{op}\otimes \mathcal{A}\rightarrow \mathcal{V}$.
Suppose, in addition, that $\mathcal{V}$ is cocomplete.  
$\mathcal{V}\text{-}{\bf Mod}(\mathcal{A},\mathcal{A})$
is a biclosed monoidal $\mathcal{V}$-category, 
with monoidal product 
$$\phi \circ \psi (a,a')=\overset{x\in Ob(\mathcal{A})} \int \phi(a,x)\otimes \psi(x,a')$$
unit $\mathcal{A}:\mathcal{A}^{op}\otimes\mathcal{A}\rightarrow \mathcal{V}$
and right internal hom 
$$[\phi,\psi]_{r}(a,a')=\underset{x\in Ob(\mathcal{A})} \int \psi(x,a')^{ \phi(x,a)}$$

Let $D:Set\rightarrow \mathcal{V}\text{-}{\bf Cat}$ be the 
discrete $\mathcal{V}$-category functor. We denote 
$\mathcal{V}\text{-}{\bf Mod}(DS,DS)$ by 
$\mathcal{V}$\text{-}{\bf Graph}$(S)$--
this is just the functor category $\mathcal{V}^{S\times S}$.
 The objects of $\mathcal{V}$\text{-}{\bf Graph}$(S)$ are called
$\mathcal{V}$-graphs with fixed set of objects $S$. For 
$\mathcal{X},\mathcal{Y}\in \mathcal{V}\text{-}{\bf Graph}(S)$, 
one now has $$\mathcal{X}\circ \mathcal{Y}(a,b)=
\underset{z\in S}\coprod \mathcal{X}(a,z)\otimes
\mathcal{Y}(z,b)$$ the unit is
$$\mathcal{I}_{S}(a,b) =
  \begin{cases}
    I, & \text{if }  a=b\\
    \emptyset, & \text{otherwise}
  \end{cases} $$ 
and the right internal hom is
$$[\mathcal{X},\mathcal{Y}]_{r}(a,b)=\underset{x\in S} 
\prod \mathcal{Y}(x,b)^{\mathcal{X}(x,a)}$$
There is an adjunction
$$\delta:\mathcal{V}\rightleftarrows 
\mathcal{V}\text{-}{\bf Graph}(S):\underline{\underline{()}}$$
where $$\delta_{X}(a,b)=
 \begin{cases}
    X, & \text{if }  a=b\\
    \emptyset, & \text{otherwise}
  \end{cases} $$
and $\underline{\underline{\mathcal{X}}}=
\underset{a\in Ob(\mathcal{A})}\prod \mathcal{X}(a,a)$. 
The functor $\delta$ is strong monoidal, therefore 
$\underline{\underline{()}}$ is monoidal for formal reasons.

The category $\mathcal{V}\text{-}{\bf Cat}(S)$ is precisely 
the category of monoids in $\mathcal{V}\text{-}{\bf Graph}(S)$ 
with respect to $\circ$, and 
$\mathcal{V}\text{-}{\bf Mod}(\mathcal{A},\mathcal{A})$
is precisely the category of $(\mathcal{A},\mathcal{A})$-bimodules in
$(\mathcal{V}\text{-}{\bf Graph}(Ob(\mathcal{A})),\circ)$ (in the sense
of categorical algebra). Thus, there is a free-forgetful adjunction
$$F=\mathcal{A}\circ-\circ \mathcal{A}:\mathcal{V}\text{-}{\bf Graph}(Ob(\mathcal{A}))
\rightleftarrows \mathcal{V}\text{-}{\bf Mod}(\mathcal{A},\mathcal{A}):U$$
One has $F\mathcal{X}\circ F\mathcal{Y}\cong    
F(\mathcal{X}\circ \mathcal{A}\circ \mathcal{Y})$
for every $\mathcal{V}$-graph $\mathcal{X}$ and every 
$\phi \in \mathcal{V}\text{-}{\bf Mod}(\mathcal{A},\mathcal{A})$,
which implies that $$U[F\mathcal{X},\phi]_{r}\cong 
[\mathcal{A}\circ \mathcal{X},U\phi]_{r}$$ Moreover,
$U(\phi \circ \psi)$ is the coequaliser of the (reflexive) pair
\[
   \xymatrix{
U(\phi)\circ \mathcal{A}\circ U(\psi)
\ar @<0pt> [r] \ar @<6pt> [r] & U(\phi)\circ U(\psi)\\
  }
  \]
so that 

$(a)$ the forgetful functor $U$ is monoidal, and

$(b)$ $[\phi,\psi]_{r}$ is the equaliser of the pair
\[
   \xymatrix{
[U(\phi),U(\psi)]_{r} \ar @<0pt> [r] \ar @<6pt> [r] & 
[U(\phi)\circ \mathcal{A},U(\psi)]_{r}\\
 }
  \]
There is a functor $$C:\mathcal{V}\text{-}{\bf Mod}
(\mathcal{A},\mathcal{A})\times \Delta_{+}^{op}\rightarrow 
\mathcal{V}\text{-}{\bf Mod}(\mathcal{A},\mathcal{A})$$
given by $$C(\phi,n)=U\phi \circ \mathcal{A}^{\circ n}
\cong \phi \circ \mathcal{A}^{\circ (n+1)}$$
The functor $C(\phi,-)$ is (strong) monoidal. Let us 
denote by $\bigstar$ the Day convolution product on 
$\mathcal{V}\text{-}{\bf Mod}(\mathcal{A},\mathcal{A})^{\Delta_{+}^{op}}$;
its unit object is $\underset{\Delta_{+}(-,0)}\sqcup \mathcal{A}$. One has
$$(\phi\bigstar \psi)(n)\cong \underset{i+j=n}\coprod \phi(i)\circ \psi(j)$$
Setting $C(\phi)(n)=C(\phi,n)$ defines a monoidal functor
$$C: \mathcal{V}\text{-}{\bf Mod}(\mathcal{A},\mathcal{A})
\rightarrow \mathcal{V}\text{-}{\bf Mod}(\mathcal{A},\mathcal{A})
^{\Delta_{+}^{op}}$$ In particular, $C(\mathcal{A})$ is a monoid in 
$\mathcal{V}\text{-}{\bf Mod}(\mathcal{A},\mathcal{A})^{\Delta_{+}^{op}}$. 
>From the last paragraph of section 3 it follows that, in the adjoint pair 
$$-\bigstar_{\Delta_{+}}C(\mathcal{A}):\mathcal{V}\text{-}
{\bf Mod}(\mathcal{A},\mathcal{A})^{\Delta_{+}}\rightleftarrows   
\mathcal{V}\text{-}{\bf Mod}(\mathcal{A},\mathcal{A}):[C(\mathcal{A}),-]_{r}$$
the left adjoint is monoidal. It can be readily seen 
that $C(\mathcal{A})$ is a non-counital comonoid in 
$\mathcal{V}\text{-}{\bf Mod}
(\mathcal{A},\mathcal{A})^{\Delta_{+}^{op}}$.
\begin{lemma}
The functor $[C(\mathcal{A}),-]_{r}$ is monoidal.
\end{lemma}
\begin{proof}
The THC-transpose of the natural map 
$(\underset{\Delta_{+}(0,-)}\sqcup \mathcal{A})
\bigstar_{\Delta_{+}}C(\mathcal{A})\rightarrow \mathcal{A}$ 
is the unit map. We shall construct a map
$$F_{\phi,\psi}:[C(\mathcal{A}),\phi]_{r}\bigstar [C(\mathcal{A}),\psi]_{r}
\rightarrow [C(\mathcal{A})\bigstar C(\mathcal{A}),\phi \circ \psi]_{r}$$
>From the previous considerations we have $U([C(\mathcal{A}),\phi]_{r}(n))=
[\mathcal{A}^{\circ n},U\phi]_{r}$.
Using this, one can define a natural cup product
$$\smile:[\mathcal{A}^{\circ m},U\phi]_{r}\circ 
[\mathcal{A}^{\circ n},U\psi]_{r}\rightarrow 
[\mathcal{A}^{\circ (m+n)},U(\phi \circ \psi)]_{r}$$
which is compatible with the actions of $\mathcal{A}$ and is 
suitably associative and unital,
so that it induces a suitably associative and unital map
$$[C(\mathcal{A})(m),\phi]_{r}\circ[C(\mathcal{A})(n),\psi]_{r}
\rightarrow [C(\mathcal{A})(m+n),\phi \circ \psi]_{r}$$
This map induces the map $F_{\phi,\psi}$. It follows that 
$[C(\mathcal{A}),-]_{r}$ is associative and unital. 
\end{proof}
\begin{notation}
{\rm We denote by $Y_{+}$ the composite 
$$\mathcal{V}\text{-}{\bf Mod}(\mathcal{A},\mathcal{A})\overset{[C(\mathcal{A}),-]_{r}}
\longrightarrow \mathcal{V}\text{-}{\bf Mod}(\mathcal{A},\mathcal{A})^{\Delta_{+}}
\overset{U^{\Delta_{+}}}\longrightarrow 
\mathcal{V}\text{-}{\bf Graph}(Ob(\mathcal{A}))^{\Delta_{+}}
\overset{\underline{\underline{()}}^{\Delta_{+}}}
\longrightarrow \mathcal{V}^{\Delta_{+}}$$ and by $Y$ 
the functor $i^{\ast}Y_{+}$, where $i^{\ast}:\mathcal{V}^{\Delta_{+}}
\rightarrow \mathcal{V}^{\Delta}$ is the restriction functor.}
\end{notation}
It follows from lemma 4.1, the previous considerations and 2.2
that $Y$ is monoidal. The functor $Y$ is a familiar one. For $\phi \in 
\mathcal{V}\text{-}{\bf Mod}(\mathcal{A},\mathcal{A})$, the cosimplicial
object $Y(\phi)^{\bullet}$ in $\mathcal{V}$ is given by
$$Y(\phi)^{n}=
\begin{cases}
\underset{a\in Ob(\mathcal{A})}\prod \phi(a,a) , & \text{if }  n=0\\
\underset{a_{0},...,a_{n}\in Ob(\mathcal{A})}\prod \phi(a_{0},a_{n})
^{\underline{\mathcal{A}}(a_{0},...,a_{n})}, & \text{if } n\geq 1,
\end{cases} $$
where $\underline{\mathcal{A}}(a_{0},...,a_{n})=\mathcal{A}(a_{0},a_{1})
\otimes...\otimes \mathcal{A}(a_{n-1},a_{n})$. Here are some examples
of coface and codegeneracy maps, for a full description see \cite[pages 6 and 7]{CP}.
$d^{0}:Y(\phi)^{n}\rightarrow Y(\phi)^{n+1}$ is obtained from the diagram
\[
   \xymatrix{
\underset{a_{0},...,a_{n}\in Ob(\mathcal{A})}\prod 
\phi(a_{0},a_{n})^{\underline{\mathcal{A}}(a_{0},...,a_{n})}
\ar @{-->} [rr] \ar[d]^{pr_{b_{1},...,b_{n+1}}} & &
\underset{b_{0},...,b_{n+1}\in Ob(\mathcal{A})}
\prod \phi(b_{0},b_{n+1})^{\underline{\mathcal{A}}(b_{0},...,b_{n+1})}
\ar[d]^{pr_{b_{0},...,b_{n+1}}}\\
\phi(b_{1},b_{n+1})^{\underline{\mathcal{A}}(b_{1},...,b_{n+1})} \ar[rr] & &
\phi(b_{0},b_{n+1})^{\underline{\mathcal{A}}(b_{0},...,b_{n+1})},\\
}
   \]
where the bottom horizontal map is the adjoint transpose of
$$\phi(b_{1},b_{n+1})^{\underline{\mathcal{A}}(b_{1},...,b_{n+1})} 
\otimes \mathcal{A}(b_{0},b_{1}) \otimes
\underline{\mathcal{A}}(b_{1},...,b_{n+1})\rightarrow 
\mathcal{A}(b_{0},b_{1}) \otimes \phi(b_{1},b_{n+1})
\rightarrow$$
$$(\mathcal{A}^{op}\otimes \mathcal{A})((b_{1},b_{n+1}),(b_{0},b_{n+1}))
\otimes \phi(b_{1},b_{n+1})\rightarrow \phi(b_{0},b_{n+1})$$
whereas $d^{n+1}$ is obtained from the diagram
\[
   \xymatrix{
\underset{a_{0},...,a_{n}\in Ob(\mathcal{A})}\prod 
\phi(a_{0},a_{n})^{\underline{\mathcal{A}}(a_{0},...,a_{n})}
\ar @{-->} [rr] \ar[d]^{pr_{b_{0},...,b_{n}}} & &
\underset{b_{0},...,b_{n+1}\in Ob(\mathcal{A})}\prod 
\phi(b_{0},b_{n+1})^{\underline{\mathcal{A}}(b_{0},...,b_{n+1})}
\ar[d]^{pr_{b_{0},...,b_{n+1}}}\\
\phi(b_{0},b_{n})^{\underline{\mathcal{A}}(b_{0},...,b_{n})} \ar[rr] & &
\phi(b_{0},b_{n+1})^{\underline{\mathcal{A}}(b_{0},...,b_{n+1})}\\
}
   \]
where the bottom horizontal map is the adjoint transpose of
$$\phi(b_{0},b_{n})^{\underline{\mathcal{A}}(b_{0},...,b_{n})} \otimes
\underline{\mathcal{A}}(b_{0},...,b_{n}) \otimes 
\mathcal{A}(b_{n},b_{n+1})\rightarrow \phi(b_{0},b_{n}) \otimes
\mathcal{A}(b_{n},b_{n+1}) \rightarrow$$
$$(\mathcal{A}^{op}\otimes 
\mathcal{A})((b_{0},b_{n}),(b_{0},b_{n+1}))\otimes 
\phi(b_{0},b_{n})\rightarrow \phi(b_{0},b_{n+1})$$
Similarly, the codegeneracy 
$s^{i}:Y(\phi)^{n+1}\rightarrow Y(\phi)^{n}$ is obtained from the diagram
\[
   \xymatrix{
\underset{a_{0},...,a_{n+1}\in Ob(\mathcal{A})}\prod 
\phi(a_{0},a_{n})^{\underline{\mathcal{A}}(a_{0},...,a_{n+1})}
\ar @{-->} [rr] \ar[d]^{pr_{b_{0},...,b_{i-1},b_{i},b_{i},...,b_{n}}} & &
\underset{b_{0},...,b_{n}\in Ob(\mathcal{A})}\prod 
\phi(b_{0},b_{n})^{\underline{\mathcal{A}}(b_{0},...,b_{n})}
\ar[d]^{pr_{b_{0},...,b_{n}}}\\
\phi(b_{0},b_{n})^{\underline{\mathcal{A}}(b_{0},...,b_{i})\otimes
\mathcal{A}(b_{i},b_{i}) \otimes
\underline{\mathcal{A}}(b_{i},...,b_{n})} \ar[rr]^{insert \ id_{b_{i}}} & &
\phi(b_{0},b_{n})^{\underline{\mathcal{A}}(b_{0},...,b_{n})}\\
}
 \]
\begin{calculation}
{\rm In 5.1 we shall need an understanding of the category 
$(\mathcal{V}\text{-}{\bf Mod}(\mathcal{A},\mathcal{A}))\text{-}{\bf Cat}$.

\emph{Step 1.} Let $(\mathcal{E},\otimes,I)$ be an arbitrary monoidal category 
having sufficient colimits and with monoidal product preserving the
existent colimits in each variable separately. We denote by $Mon(\mathcal{E})$
the category of monoids in $\mathcal{E}$. Let $A\in Mon(\mathcal{E})$. 
Let $_{A}Mod_{A}$ be the category of $(A,A)$-bimodules in $\mathcal{E}$.
The categories $Mon(_{A}Mod_{A})$ and $(A\downarrow Mon(\mathcal{E}))$
are isomorphic as categories above $Mon(\mathcal{E})$:
\[
   \xymatrix{
Mon(_{A}Mod_{A}) \ar[rr]^{\cong} \ar[dr] & & 
(A\downarrow Mon(\mathcal{E})) \ar[dl]\\
& Mon(\mathcal{E})\\
}
  \]
\emph{Step 2.} Let $(\mathcal{E},\otimes,I)$ be as in step 1 and 
let $S$ be a set. At the beginning of this section we defined a strong 
monoidal functor $\delta:\mathcal{E}\rightarrow \mathcal{E}\text{-}{\bf Graph}(S)$. 
Let $A$ be a monoid in $\mathcal{E}$. The categories 
$_{A}Mod_{A}\text{-}{\bf Cat}(S)$ and
$Mon(_{\delta_{A}}Mod_{\delta_{A}})$
are isomorphic. Therefore, by step 1 the categories 
$_{A}Mod_{A}\text{-}{\bf Cat}(S)$ and 
$(\delta_{A}\downarrow \mathcal{E}\text{-}{\bf Cat}(S))$
are isomorphic as categories above $\mathcal{E}\text{-}{\bf Cat}(S)$:
\[
   \xymatrix{
_{A}Mod_{A}\text{-}{\bf Cat}(S) \ar[rr]^{\cong} \ar[dr] & & 
(\delta_{A}\downarrow \mathcal{E}\text{-}{\bf Cat}(S)) \ar[dl]\\
& \mathcal{E}\text{-}{\bf Cat}(S)\\
}
  \]
\emph{Step 3.} Let $(\mathcal{V},\otimes,I)$ be a closed category 
having sufficient colimits. Let $\mathcal{A}$ be a $\mathcal{V}$-category
and $S$ a set. By step 2 and previous considerations the categories
$(\mathcal{V}\text{-}{\bf Mod}(\mathcal{A},\mathcal{A}))\text{-}{\bf Cat}(S)$
and $(\delta_{\mathcal{A}}\downarrow (\mathcal{V}\text{-}{\bf Graph}
(Ob(\mathcal{A})))\text{-}{\bf Cat}(S))$ are isomorphic as categories 
above $(\mathcal{V}\text{-}{\bf Graph}(Ob(\mathcal{A})))\text{-}{\bf Cat}(S)$:
\[
   \xymatrix{
(\mathcal{V}\text{-}{\bf Mod}(\mathcal{A},\mathcal{A}))\text{-}{\bf Cat}(S)
\ar[rr]^{\cong} \ar[dr] & & 
(\delta_{\mathcal{A}}\downarrow (\mathcal{V}\text{-}{\bf Graph}
(Ob(\mathcal{A})))\text{-}{\bf Cat}(S)) \ar[dl]\\
& (\mathcal{V}\text{-}{\bf Graph}(Ob(\mathcal{A})))\text{-}{\bf Cat}(S)\\
}
  \]
The right-hand corner of the previous diagram implies that
to give an object of $(\mathcal{V}\text{-}{\bf Mod}
(\mathcal{A},\mathcal{A}))\text{-}{\bf Cat}$ amounts to 
a set $S$, a $\mathcal{V}$-category $\mathcal{Z}$
with object set $S\times Ob(\mathcal{A})$ and, for all
$x\in S$, $\mathcal{V}$-functors $u_{x}:\mathcal{A}
\rightarrow \mathcal{Z}$ with object maps $u_{x}(a)=(x,a)$.}
\end{calculation}

\section{$\mathcal{V}$-coherent transformations with 
respect to a coaugmented cosimplicial comonoid and 
the Gray tensor product with respect to a 
cosimplicial comonoid}

Throughout this section $(\mathcal{V},\otimes,I)$ is a complete 
and cocomplete closed category. We write $Y^{X}$ for the 
internal hom of two objects $X$, $Y$ of $\mathcal{V}$.
Given two small $\mathcal{V}$-categories $\mathcal{A}$, 
$\mathcal{B}$ and a comonoid $C$ in $\mathcal{V}^{\Delta_{+}}$, 
we shall construct a $\mathcal{V}$-category 
$Coh^{C}(\mathcal{A},\mathcal{B})$ whose objects are the 
$\mathcal{V}$-functors from $\mathcal{A}$ to $\mathcal{B}$.
Of particular importance will be the comonoids $C$ of the form 
$i_{!}C^{\bullet}$ (2.2).

Let  $C$ be a comonoid in $\mathcal{V}^{\Delta_{+}}$
and $\mathcal{A}$ a $\mathcal{V}$-category. We define
$Coh^{C}(\mathcal{A},-):\mathcal{V}\text{-}{\bf Cat}\rightarrow
\mathcal{V}\text{-}{\bf Cat}$ as the composite of the three canonical functors
$$\mathcal{V}\text{-}{\bf Cat}\underset{\mathcal{V}^{\Delta_{+}}(C,-)}
\longleftarrow \mathcal{V}^{\Delta_{+}}\text{-}{\bf Cat}\underset{Y_{+}}
\longleftarrow (\mathcal{V}\text{-}{\bf Mod}(\mathcal{A},\mathcal{A}))\text{-}{\bf Cat}
\underset{\langle\mathcal{A},-\rangle}\longleftarrow \mathcal{V}\text{-}{\bf Cat}$$
The arrows point to the right in order to emphasize that each of the three
functors is a right adjoint, as we shall later show. We need to construct
$\langle\mathcal{A},-\rangle$. If $\mathcal{B}$ is a $\mathcal{V}$-category, 
the objects of $\langle\mathcal{A},\mathcal{B}\rangle$ 
are the $\mathcal{V}$-functors $\mathcal{A}\rightarrow \mathcal{B}$ and
$$\langle\mathcal{A},\mathcal{B}\rangle(f,g)(a,a')=
\mathcal{B}(f^{op}\otimes g)(a,a')=\mathcal{B}(fa,ga')$$
The composition maps $\langle\mathcal{A},\mathcal{B}\rangle(f,g)\circ
\langle\mathcal{A},\mathcal{B}\rangle(g,h)\rightarrow 
\langle\mathcal{A},\mathcal{B}\rangle(f,h)$
are induced by the composition maps of $\mathcal{B}$. Varying $C$ and 
$\mathcal{A}$ we obtain a functor
$$Coh^{-}(-,-):(Comon(\mathcal{V}^{\Delta_{+}})\times 
\mathcal{V}\text{-}{\bf Cat})^{op}
\times \mathcal{V}\text{-}{\bf Cat}\rightarrow \mathcal{V}\text{-}{\bf Cat}$$
$$Coh^{C}(\mathcal{A},\mathcal{B})=
\mathcal{V}^{\Delta_{+}}(C,Y_{+}\langle\mathcal{A},\mathcal{B}\rangle)$$
We call $Coh^{C}(\mathcal{A},\mathcal{B})(f,g)$ the {\bf object
of $\mathcal{V}$-coherent transformations from $f$ to $g$ with respect to $C$}.
Let $C^{\bullet}$ be a comonoid in $\mathcal{V}^{\Delta}$.
We call $Coh^{i_{!}C^{\bullet}}(\mathcal{A},\mathcal{B})(f,g)$ 
the {\bf object of $\mathcal{V}$-coherent transformations from $f$ to 
$g$ with respect to $C^{\bullet}$}.
We shall write $Coh^{C^{\bullet}}(\mathcal{A},\mathcal{B})$ instead of
$Coh^{i_{!}C^{\bullet}}(\mathcal{A},\mathcal{B})$.
Since the adjunction $$-\ast i_{!}C^{\bullet}:\mathcal{V}\rightleftarrows
\mathcal{V}^{\Delta_{+}}:\mathcal{V}^{\Delta_{+}}(i_{!}C^{\bullet},-)$$
splits as
\[
   \xymatrix{
\mathcal{V} \ar @<2pt> [r]^{-\ast C^{\bullet}} & \mathcal{V}^{\Delta} 
\ar @<2pt> [l]^{\mathcal{V}^{\Delta}(C^{\bullet},-)}
\ar @<2pt> [r]^{i_{!}} & \mathcal{V}^{\Delta_{+}} 
\ar @<2pt> [l]^{i^{\ast}}\\
  }
  \]
it follows that $$Coh^{C^{\bullet}}(\mathcal{A},\mathcal{B})=
\mathcal{V}^{\Delta}(C^{\bullet},Y\langle\mathcal{A},\mathcal{B}\rangle)$$
The previous formula was considered by Cordier and Porter \cite[Definition 3.1]{CP}
in the case when $\mathcal{V}$ is the category of simplicial sets and
the cosimplicial simplicial set $\Delta$ is in the place of $C^{\bullet}$. See also the
references therein. We have borrowed the notation $Coh$ from them. However, 
the formalism leading to this formula is not present in \cite{CP}. Also, in the case
of simplicial sets it is known that $\Delta$ is not a comonoid with respect to
$\bigstar$ and that $Coh(\mathcal{A},\mathcal{B})$ cannot be naturally made
into a simplicial category; see \cite[page 28]{CP} for a discussion of the latter fact.
\begin{examples}
$(a)$ {\rm Given a $\mathcal{V}$-category $\mathcal{A}$ and a 
comonoid $C$ in $\mathcal{V}$, let us denote by 
$\mathcal{A}^{C}$ the $\mathcal{V}$-category having the same 
objects as $\mathcal{A}$ and having the $\mathcal{V}$-homs 
$\mathcal{A}^{C}(a,a')=\mathcal{A}(a,a')^{C}$. Then one has
$Coh^{C}(\mathcal{I},\mathcal{A})\cong \mathcal{A}^{C(1)}$ and 
$Coh^{C^{\bullet}}(\mathcal{I},\mathcal{A})\cong \mathcal{A}^{C^{0}}$.}

$(b)$ {\rm $Coh^{cstI}(\mathcal{A},\mathcal{B})$ coincides with the 
internal hom $\mathcal{V}Nat(\mathcal{A},\mathcal{B})$ of the 
standard closed category structure on $\mathcal{V}\text{-}{\bf Cat}$, since 
$Coh^{cstI}(\mathcal{A},\mathcal{B})(f,g)$
is the object of $\mathcal{V}$-natural transformations from $f$ to $g$.}

$(c)$ {\rm $Coh^{sk(I)}(\mathcal{A},\mathcal{B})(f,g)=
\underset{a\in Ob(\mathcal{A})}\prod \mathcal{B}(fa,ga)$. More generally, 
for a comonoid $C$ in $\mathcal{V}$, $Coh^{sk(C)}(\mathcal{A},\mathcal{B})=
Coh^{sk(I)}(\mathcal{A},\mathcal{B})^{C}$.}
\end{examples}
\begin{remark} {\rm There are variants of $Coh^{C^{\bullet}}(\mathcal{A},\mathcal{B})$,
if one is willing to replace $\mathcal{V}^{\Delta}$ with $\mathcal{V}^{\Delta(1)}$
or $\mathcal{V}^{\Delta(1)_{mon}}$. For example, in the case of
$\mathcal{V}^{\Delta(1)}$ one obtains a functor
$$Coh^{-}_{\Delta(1)}(-,-): (Comon(\mathcal{V}^{\Delta(1)})\times
\mathcal{V}\text{-}{\bf Cat})^{op}\times \mathcal{V}\text{-}{\bf Cat}
\rightarrow \mathcal{V}\text{-}{\bf Cat}$$
In the case of $\mathcal{V}^{\Delta(1)_{mon}}$
 one obtains a $\mathcal{V}$-category without unit.}
\end{remark}
\subsection{The main TH-situation}
We shall now show that $Coh^{-}(-,-)$ is part of a TH-situation
$$-\boxtimes_{-}-:\mathcal{V}{\bf Cat}\times
(Comon(\mathcal{V}^{\Delta_{+}})\times \mathcal{V}\text{-}{\bf Cat})
\rightarrow \mathcal{V}\text{-}{\bf Cat}$$
$$Coh^{-}(-,-):(Comon(\mathcal{V}^{\Delta_{+}})\times
\mathcal{V}\text{-}{\bf Cat})^{op}\times 
\mathcal{V}\text{-}{\bf Cat}\rightarrow \mathcal{V}\text{-}{\bf Cat}$$
To construct this TH-situation it suffices to construct, for every
comonoid $C$ in $\mathcal{V}^{\Delta_{+}}$, a TH-situation
$$-\boxtimes_{C}-:\mathcal{V}\text{-}{\bf Cat}\times\mathcal{V}\text{-}{\bf Cat}
\rightarrow \mathcal{V}\text{-}{\bf Cat}$$
$$Coh^{C}(-,-):\mathcal{V}\text{-}{\bf Cat}^{op}\times 
\mathcal{V}\text{-}{\bf Cat}\rightarrow 
\mathcal{V}\text{-}{\bf Cat}$$ 
such that the adjunction isomorphisms in this 
TH-situation are natural in $C$. In turn, to construct the latter 
TH-situation it suffices to construct, for every 
$\mathcal{V}$-category $\mathcal{A}$, a left adjoint
$-\boxtimes_{C}\mathcal{A}$ to $Coh^{C}(\mathcal{A},-)$
such that the adjunction isomorphisms are natural in $\mathcal{A}$.
We show that each of the three functors which make up 
$Coh^{C}(\mathcal{A},-)$ has a left adjoint, thus
$-\boxtimes_{C}\mathcal{A}$ will be by definition the 
composite of these left adjoints.

$\bullet$ The left adjoints to $\mathcal{V}^{\Delta_{+}}(C,-)$
and $Y_{+}$ are constructed using the general
\begin{fact} Let $\mathcal{E}_{1}$ and $\mathcal{E}_{2}$ be two 
cocomplete monoidal categories with monoidal products cocontinuous
in each variable separately. We denote by 
$\mathcal{E}_{i}\text{-}{\bf Graph}$ the category of 
small $\mathcal{E}_{i}$-graphs. Let 
$\mathcal{F}_{i}:\mathcal{E}_{i}\text{-}{\bf Graph}\rightleftarrows
\mathcal{E}_{i}\text{-}{\bf Cat}:\mathcal{U}_{i}$
be the free-forgetful adjunction \cite{Wo}, $i\in \{1,2\}$. Let
$F:\mathcal{E}_{1}\rightleftarrows \mathcal{E}_{2}:G$
be an adjoint pair with $G$ monoidal. The functor 
$G: \mathcal{E}_{2}\text{-}{\bf Cat} \rightarrow  
\mathcal{E}_{1}\text{-}{\bf Cat}$ has a left adjoint $F'$ 
constructed in such a way that $\mathcal{F}_{2}F
\cong F'\mathcal{F}_{1}$. In particular, $F'$ preserves 
the unit object. $F'$ is constructed (fibrewise) as follows.
For $\mathcal{A}\in \mathcal{E}_{1}\text{-}{\bf Cat}$,
$F'\mathcal{A}$ is the coequaliser of the reflexive pair
\[
   \xymatrix{
\mathcal{F}_{2}F\mathcal{U}_{1}\mathcal{F}_{1}\mathcal{U}_{1}\mathcal{A}
\ar @<-2pt> [r] \ar @<2pt> [r] & 
\mathcal{F}_{2}F\mathcal{U}_{1}\mathcal{A}\\
 }
  \]
One arrow is obtained by applying $\mathcal{F}_{2}F\mathcal{U}_{1}$
to the counit $\mathcal{F}_{1}\mathcal{U}_{1}\mathcal{A}\rightarrow \mathcal{A}$.
The other one is obtained by substituting 
$\mathcal{X}=\mathcal{U}_{1}\mathcal{A}$
in the adjoint transpose of the natural map 
$F\mathcal{U}_{1}\mathcal{F}_{1}\mathcal{X}\rightarrow 
\mathcal{U}_{2}\mathcal{F}_{2}F\mathcal{X}$, 
$\mathcal{X}\in \mathcal{E}_{1}\text{-}{\bf Graph}$.
\end{fact}
$\bullet$ The left adjoint $-\lozenge\mathcal{A}$ to $\langle\mathcal{A},-\rangle$
is defined as follows. $Ob(\mathcal{C}\lozenge\mathcal{A})=
Ob(\mathcal{C})\times Ob(\mathcal{A})$ and 
$\mathcal{C}\lozenge\mathcal{A}((c,a),(c',a'))=\mathcal{C}(c,c')(a,a')$.
To see that $\mathcal{C}\lozenge \mathcal{A}$ is well-defined
and indeed a left adjoint one uses calculation 4.3.

This finishes the construction of $-\boxtimes_{C}\mathcal{A}$.
The adjunction isomorphism is clearly natural in $\mathcal{A}$, 
hence we obtain a TH-situation
$$(-\bigstar_{\Delta_{+}}C(-))(F^{\Delta_{+}})'\delta^{\Delta_{+}}(-)
\lozenge-: \mathcal{V}^{\Delta_{+}}\text{-}{\bf Cat}
\times\mathcal{V}\text{-}{\bf Cat}\rightarrow 
\mathcal{V}\text{-}{\bf Cat}$$
$$Y_{+}\langle-,-\rangle:\mathcal{V}\text{-}{\bf Cat}^{op}
\times \mathcal{V}\text{-}{\bf Cat}
\rightarrow \mathcal{V}^{\Delta_{+}}\text{-}{\bf Cat}$$
Using the adjunction $-\ast C:\mathcal{V}\rightleftarrows
\mathcal{V}^{\Delta_{+}}:\mathcal{V}^{\Delta_{+}}(C,-)$, 
where $C$ is a comonoid in $\mathcal{V}^{\Delta_{+}}$, 
and the fact that TH-situations can be changed along 
adjoint functors, we obtain a TH-situation
$$-\boxtimes_{C}-:\mathcal{V}\text{-}{\bf Cat}\times
\mathcal{V}\text{-}{\bf Cat}\rightarrow \mathcal{V}\text{-}{\bf Cat}$$
$$Coh^{C}(-,-):\mathcal{V}\text{-}{\bf Cat}^{op}\times 
\mathcal{V}\text{-}{\bf Cat}\rightarrow 
\mathcal{V}\text{-}{\bf Cat}$$ 
It is clear that the adjunction isomorphism in this TH-situation
is natural in $C$, therefore we obtain the desired TH-situation.
Using again the fact that TH-situations can be changed along
adjoint functors, we obtain a TH-situation
$$-\boxtimes_{-}-:\mathcal{V}{\bf Cat}\times
(Comon(\mathcal{V}^{\Delta})\times \mathcal{V}\text{-}{\bf Cat})
\rightarrow \mathcal{V}\text{-}{\bf Cat}$$
$$Coh^{-}(-,-):(Comon(\mathcal{V}^{\Delta})\times
\mathcal{V}\text{-}{\bf Cat})^{op}\times 
\mathcal{V}\text{-}{\bf Cat}\rightarrow \mathcal{V}\text{-}{\bf Cat}$$
For a comonoid $C^{\bullet}$ in $\mathcal{V}^{\Delta}$
one has $\mathcal{A}\boxtimes_{C^{\bullet}}\mathcal{B}=
\mathcal{A}\boxtimes_{i_{!}C^{\bullet}}\mathcal{B}$.
We call  $\mathcal{A}\boxtimes_{C^{\bullet}}\mathcal{B}$
the {\bf Gray tensor product of $\mathcal{A}$ and 
$\mathcal{B}$ with respect to} $C^{\bullet}$.
We do not claim that $\boxtimes_{C^{\bullet}}$ 
is a monoidal product on $\mathcal{V}$\text{-}{\bf Cat},
but see proposition 5.9. The naming will be justified in 5.3.
\begin{notation}
For an object $X$ of an arbitrary monoidal category 
$\mathcal{E}$  with unit $I$ and having an initial object, 
we denote by $2_{X}$ the $\mathcal{E}$-category with 
two objects $0$ and $1$ and with $2_{X}(0,0)=2_{X}(1,1)=I$, 
$2_{X}(1,0)=\emptyset$ and $2_{X}(0,1)=X$.
\end{notation}
For example, in the setting of 5.3 one has $2_{FX}\cong F'(2_{X})$
for every $X\in \mathcal{E}_{1}$. If $\phi \in 
\mathcal{V}\text{-}{\bf Mod}(\mathcal{A},\mathcal{A})$
we represent $2_{\phi}\lozenge\mathcal{A}$ as
\[
   \xymatrix{
(0,a') \ar@{-} [rr]^{\phi(a',a')} & & (1,a')\\
(0,a) \ar@{-} [rr]^{\phi(a,a)} \ar@{-} [u]^{\mathcal{A}(a,a')} 
\ar@{-} [urr]^{\phi(a,a')} & & (1,a) \ar @{-} [u]_{\mathcal{A}(a,a')}\\
 }
  \]
\begin{examples}
$(a)$ {\rm $\mathcal{I}\boxtimes_{C}\mathcal{A}\cong\mathcal{A}$ and 
$-\boxtimes_{C}\mathcal{I}\cong(-\otimes C(1))'$.

$(b)$ By example 5.1$(c)$, $\mathcal{A}\boxtimes_{sk(C)}\mathcal{B}\cong
(-\otimes C)'(\mathcal{A})\boxtimes_{sk(I)}\mathcal{B}$ for an 
arbitrary comonoid $C$ in $\mathcal{V}$.

$(c)$ $2_{X}\boxtimes_{C}\mathcal{A}\cong 
2_{(\mathcal{A}\circ \delta_{X\ast C}\circ \mathcal{A})
\bigstar_{\Delta_{+}}C(\mathcal{A})}\lozenge\mathcal{A}$.}
\end{examples}
\begin{remark}
{\rm Let $C$ be a comonoid in $\mathcal{V}^{\Delta_{+}}$ 
with $C(1)\neq I$. Example 5.5$(a)$ shows that there is no 
right closed category structure on $\mathcal{V}\text{-}{\bf Cat}$ 
with unit $\mathcal{I}$ and right internal 
hom $Coh^{C}(-,-)$.}
\end{remark}
The next result, whose proof is left to the reader, 
is inspired by \cite[Section 7]{CP}.
\begin{lemma}
There are two natural maps
$$Y_{+}\langle\mathcal{A},\mathcal{B}\rangle(f,g)\bigstar
Y_{+}\langle\mathcal{B},\mathcal{C}\rangle(k,l) \rightarrow 
Y_{+}\langle\mathcal{A},\mathcal{C}\rangle(kf,lg)$$
which are suitably associative and unital. Consequently,
for every comonoid $C$ in $\mathcal{V}^{\Delta_{+}}$, $Coh^{C}(-,-)
\in (\mathcal{V}\text{-}{\bf Cat})\text{-}{\bf CAT}$ in two ways.
\end{lemma}
\subsection{The case of $sk(I)$}
It is known, more or less from \cite{FKL},
that $Coh^{sk(I)}(-,-)$ is the internal hom of
a closed category structure on $\mathcal{V}\text{-}{\bf Cat}$. 
We provide below more details than in \emph{loc.cit.}.

Let $\mathcal{A}$, $\mathcal{B}$ and $\mathcal{C}$ be three 
$\mathcal{V}$-categories. A {\bf pre}-{\bf bi}-$\mathcal{V}$-{\bf functor}
$F:(\mathcal{A},\mathcal{B})\rightarrow \mathcal{C}$
consists of the following data: for all $a\in Ob(\mathcal{A})$ 
and $b\in Ob(\mathcal{B})$ there are 
$\mathcal{V}$-functors $F(a,-):\mathcal{B}\rightarrow \mathcal{C}$
and $F(-,b):\mathcal{A}\rightarrow \mathcal{C}$ such that $F(a,-)(b)=F(-,b)(a)$.
We denote by $Pre$-$bi$-$\mathcal{V}$-$Fun(\mathcal{A},
\mathcal{B};\mathcal{C})$ the set of pre-bi-$\mathcal{V}$-functors 
$(\mathcal{A},\mathcal{B})\rightarrow \mathcal{C}$. We obtain 
a functor $Pre$-$bi$-$\mathcal{V}$-$Fun(-,-;-):(\mathcal{V}\text{-}{\bf Cat}
\times \mathcal{V}\text{-}{\bf Cat})^{op}\times 
\mathcal{V}\text{-}{\bf Cat}\rightarrow Set$.
It follows from example 5.1$(c)$ that
\begin{lemma}
There is a natural bijection $$\mathcal{V}\text{-}{\bf Cat}
(\mathcal{A},Coh^{sk(I)}(\mathcal{B},\mathcal{C}))\cong
Pre\text{-}bi\text{-}\mathcal{V}\text{-}
Fun(\mathcal{A},\mathcal{B};\mathcal{C})$$
\end{lemma}
Let now $\mathcal{A}$ and $\mathcal{B}$ be 
two $\mathcal{V}$-categories and let $S=Ob(\mathcal{A})$, 
$T=Ob(\mathcal{B})$. We define $\mathcal{A}
\boxtimes_{sk(I)}\mathcal{B}$ to be the pushout of the diagram
\[
   \xymatrix{
\mathcal{I}_{S}\otimes \mathcal{I}_{T} \ar[r] 
\ar[d] & \mathcal{I}_{S}\otimes \mathcal{B} \ar[d]\\
\mathcal{A}\otimes \mathcal{I}_{T} \ar[r] & 
\mathcal{A}\boxtimes_{sk(I)}\mathcal{B}\\
}
  \]
This pushout is calculated in $\mathcal{V}\text{-}{\bf Cat}(S\times T)$. 
\begin{proposition}
The category $(\mathcal{V}\text{-}{\bf Cat},\boxtimes_{sk(I)},\mathcal{I})$ 
is a closed category, with internal hom $Coh^{sk(I)}(-,-)$.
\end{proposition}

\subsection{The relation with Gray's tensor products}

A {\bf cocategory interval} in $\mathcal{V}$ is a
cocategory object in $\mathcal{V}$ with object of 
coobjects equal to $I$. We write such a gadget as
\[
   \xymatrix{
  I^{0}=I \ar @<-4pt> [r]_{d^{1}} 
\ar @<10pt> [r]^{d^{0}} & I^{1} \ar @<-1pt> [l]_{p}
  \ar@<10pt> [r]^{i^{0}} 
\ar@<2pt> [r]^{c} \ar@<-4pt> [r]_{i^{1}} & I^{2}\\
  }
  \]
$c$ denotes the cocomposition. There is also the
obvious notion of {\bf cogroupoid interval}.

Cocategory intervals are preserved by functors which preserve 
the unit object and finite colimits. A cocategory interval as above 
is the beginning of a cosimplicial object in $\mathcal{V}$
which we denote by $I^{\bullet}$. We shall always use the 
same notation for both the data for a cocategory interval 
and the cosimplicial object it gives rise to. We shall need to be
explicit about certain coface and codegeneracy maps of $I^{\bullet}$. 
The coface maps $d^{i}:I^{1}\rightarrow I^{2}$ are $d^{0}=i^{0}, 
d^{1}=c$ and $d^{2}=i^{1}$. The coface maps 
$d^{i}:I^{2}\rightarrow I^{3}$ are depicted in the diagrams below,
in which all squares are pushouts:
\[
   \xymatrix{
I^{1} \ar[r]^{i^{0}} \ar[d]_{i^{1}} & I^{2} \ar[d]^{d^{3}}\\
I^{2} \ar[r]^{d^{0}} & I^{3}\\
}
  \]
\[
   \xymatrix{
I \ar[r]^{d^{1}} \ar[d]_{d^{0}} & I^{1} \ar[d]_{i^{0}} \ar[ddr]^{d^{0}i^{0}}\\
I^{1} \ar[r]^{i^{1}} \ar[d]_{c} & I^{2} \ar@ {-->}  [dr]_{d^{1}}\\
I^{2} \ar[rr]_{d^{3}} & & I^{3}\\
  }
  \]
\[
   \xymatrix{
I \ar[r]^{d^{0}} \ar[d]_{d^{1}} & I^{1} \ar[d]_{i^{1}} \ar[ddr]^{d^{3}i^{1}}\\
I^{1} \ar[r]^{i^{0}} \ar[d]_{c} & I^{2} \ar@ {-->}  [dr]_{d^{2}}\\
I^{2} \ar[rr]_{d^{0}} & & I^{3}\\
  }
  \]
The coface map $s^{0}:I^{2}\rightarrow I^{1}$ is the unique 
map such that $s^{0}i^{0}=1_{I^{1}}$ and $s^{0}i^{1}=d^{1}p$. 
The other coface map $s^{1}:I^{2}\rightarrow I^{1}$ is the unique 
map such that $s^{1}i^{1}=1_{I^{1}}$ and $s^{1}i^{0}=d^{0}p$. 
More details about cocategory intervals can be found in \cite{Wa}.

\begin{examples} {\rm $(a)$ The initial cocategory interval in $Set$ has
$I^{1}=\{0,1\}$ and $I^{2}=\{0,1,2\}$, with the usual coface maps.
The cocomposition $c$ is the map which omits 1. This is a cogroupoid interval.
Therefore $\mathcal{V}$ has the initial cogroupoid interval obtained 
using the functor $F:Set\rightarrow \mathcal{V}$, $F(S)=
\underset{S}\sqcup I$; this is precisely $sk(I)$ (2.2).

$(b)$ The standard cocategory interval in {\bf Cat} 
has $I^{1}=[1]$, the totally ordered set $\{0<1\}$, and $I^{2}=[2]$, 
the totally ordered set $\{0<1<2\}$. The cocomposition is the map which 
omits 1. We shall denote this cocategory interval by
$\mathbb{I}^{\bullet}$. Applying to $\mathbb{I}^{\bullet}$ the 
free groupoid functor, we obtain the standard cogroupoid interval 
$\mathbb{J}^{\bullet}$ in the category $\mathbf{Grpd}$ of
small groupoids. $\mathbb{J}^{1}$ has two objects and one arrow 
between them and $\mathbb{J}^{2}$ has three objects and one arrow 
between any two objects. Alternatively, $\mathbb{J}^{n}$ is the 
indiscrete/chaotic category on the set $I^{n}$ considered in $(a)$.
We shall view $\mathbb{J}^{\bullet}$ as living in {\bf Cat}. 
The functor $F$ from $(a)$ induces a functor
$V:{\bf Cat}\rightarrow \mathcal{V}$-{\bf Cat}, left adjoint to the 
underlying category functor. Therefore $V(\mathbb{J}^{\bullet})$ 
is a cogroupoid interval in $\mathcal{V}$-{\bf Cat}.}
\end{examples}

\begin{proposition}
$\mathbb{I}^{\bullet}$ is a comonoid in $({\bf Cat}^{\Delta},\bigstar)$.
\end{proposition}
\begin{proof}
We first show that $\mathbb{I}^{\bullet}$ is a comonoid in 
$({\bf Cat}^{\Delta(1)},\bigstar)$. We denote by $1$ the terminal category.
Since $(\mathbb{I}^{\bullet} \bigstar \mathbb{I}^{\bullet})^{0}=1\times 1$,
we choose a map $f^{0}:1\rightarrow 1\times 1$ to be the inverse 
of the left (or right) constraint of {\bf Cat} evaluated at $1$. 
The object $(\mathbb{I}^{\bullet} \bigstar \mathbb{I}^{\bullet})^{1}$
is the pushout of the diagram $$1\times \mathbb{I}^{1}
\overset{Id_{1}\times d^{0}}\longleftarrow 
1\times 1\overset{d^{1}\times Id_{I}}\longrightarrow \mathbb{I}^{1}\times 1$$
hence it is isomorphic to $\mathbb{I}^{2}$. The coface maps are 
$D^{0}=i^{0}(d^{0}\times Id_{1})$ and $D^{1}=i^{1}(Id_{1}\times d^{1})$, 
the codegeneracy $s^{0}$ being depicted in the diagram below:
\[
   \xymatrix{
1 \ar[rr]^{f^{0}} \ar @<-4pt> [dd]_{d^{1}} \ar @<10pt> [dd]^{d^{0}} 
& &1\times 1 \ar @<-4pt> [dl]_{1\times d^{0}}
\ar @<10pt> [dl]^{1\times d^{1}} \ar @<-4pt> [dr]_{d^{0}\times 1} 
\ar @<10pt> [dr]^{d^{1}\times 1}\\
& 1\times \mathbb{I}^{1} \ar[dr]^{i^{1}} \ar[ddr]_{1\times p} & & 
\mathbb{I}^{1}\times 1 \ar[dl]_{i^{0}} \ar[ddl]^{p\times 1}\\
\mathbb{I}^{1} \ar @<-1pt> [uu]_{p} & & (\mathbb{I}^{\bullet} \bigstar 
\mathbb{I}^{\bullet})^{1}=\mathbb{I}^{2} \ar[d]^{s^{0}}\\
& & 1\times 1\\
}
   \]
This forces us to set $f^{1}=c$. We have obtained a map 
$f:\mathbb{I}^{\bullet} \rightarrow \mathbb{I}^{\bullet} 
\bigstar \mathbb{I}^{\bullet}$. Next, using the description of 
$\mathbb{I}^{\bullet} \bigstar \mathbb{I}^{\bullet}$ one can see 
that $((\mathbb{I}^{\bullet} \bigstar \mathbb{I}^{\bullet})\bigstar 
\mathbb{I}^{\bullet})^{1}$ is (isomorphic to) the pushout of the diagram
\[
   \xymatrix{
\mathbb{I}^{1}\\
1 \ar[r]^{d^{1}} \ar[u]^{d^{0}} & \mathbb{I}^{1} 
\ar[r]^{i^{1}} & \mathbb{I}^{2}\\
}
  \]
hence it is $\mathbb{I}^{3}$. The coface maps are 
$D^{0}=d^{3}i^{1}d^{1}$ and $D^{1}=d^{0}i^{0}d^{0}$. 
The above pushout can also be calculated as the pushout of the diagram
\[
   \xymatrix{
\mathbb{I}^{2}\\
\mathbb{I}^{1} \ar[u]^{i^{0}}\\
1 \ar[r]^{d^{1}} \ar[u]^{d^{0}} & \mathbb{I}^{1}\\
}
  \]
which is $(\mathbb{I}^{\bullet} \bigstar (\mathbb{I}^{\bullet} 
\bigstar \mathbb{I}^{\bullet}))^{1}$. 
It follows that $(f\bigstar \mathbb{I}^{\bullet})^{1}=d^{2}$ 
and $(\mathbb{I}^{\bullet}\bigstar f)^{1}=d^{1}$, therefore 
$f$ is coassociative. The counit $\mathbb{I}^{\bullet}\rightarrow cst1$ 
is the codegeneracy map.

Using now the presentations of $\bigstar$ from 2.2 
one can compute that the category 
$(\mathbb{I}^{\bullet} \bigstar \mathbb{I}^{\bullet})^{2}$ is
\[
   \xymatrix{
a_{0,2} \ar[r] & a_{1,2} \ar[r] & a_{2,2}\\
a_{0,1} \ar[r] \ar[u] & a_{1,1} \ar[u]\\
a_{0,0} \ar[u]\\
 }
  \]
(More generally, $(\mathbb{I}^{\bullet} \bigstar \mathbb{I}^{\bullet})^{n}$ 
is the functor category $\mathbb{I}^{n^{\mathbb{I}^{1}}}$.) The cofaces 
$D^{i}:\mathbb{I}^{2}\rightarrow (\mathbb{I}^{\bullet} 
\bigstar \mathbb{I}^{\bullet})^{2}$ are given by
$D^{0}(0)=a_{1,1}$, $D^{0}(1)=a_{1,2}$, 
$D^{0}(2)=a_{2,2}$, $D^{1}(0)=a_{0,0}$, $D^{1}(1)=a_{0,2}$
$D^{1}(2)=a_{2,2}$, $D^{2}(0)=a_{0,0}$, $D^{2}(1)=a_{0,1}$, 
and $D^{2}(2)=a_{1,1}$. Then $f^{2}:=(D^{0}c,D^{2}c)$ is given 
by $f^{2}(0)=a_{0,0}$, $f^{2}(1)=a_{1,1}$ and $f^{2}(2)=a_{2,2}$, 
and so $f^{2}c=D^{1}c$. By lemma 5.12 we obtain a map 
$f^{\bullet}:\mathbb{I}^{\bullet}\rightarrow
\mathbb{I}^{\bullet} \bigstar \mathbb{I}^{\bullet}$ 
which is coassociative. The counit axiom is easy to see.
\end{proof}
Let $c^{i}:[1]\rightarrow [n]$ be the cosimplicial operator
given by $c^{i}(0)=i-1$ and $c^{i}(1)=i$, where $1\leq i\leq n$.
We write $${\bf res}:\mathcal{V}^{\Delta}\rightarrow 
\mathcal{V}^{\Delta(1)}$$ for the functor given by 
restriction along the inclusion ${\bf res}:\Delta(1)\subset \Delta$.

\begin{lemma} Let $A^{\bullet}$ be a cocategory object in $\mathcal{V}$,
$X^{\bullet}$ is a cosimplicial object in $\mathcal{V}$
and $f:{\bf res}(A^{\bullet})\rightarrow {\bf res}(X^{\bullet})$.
For $n\geq 2$ denote by $f^{n}:A^{n}\rightarrow X^{n}$
the unique morphism of $\mathcal{V}$ such that $f^{n}c^{i}
=c^{i}f^{1}$. Then

$(1)$ There is at most one $\bar{f}:A^{\bullet}\rightarrow X^{\bullet}$
such that ${\bf res}(\bar{f})=f$.

$(2)$ Such an $\bar{f}$ exists if and only if $f^{2}c=d^{1}f^{1}$,
and when this is so, $\bar{f}^{n}=f^{n}$.
\end{lemma}

Let us take $\mathcal{V}=({\bf Cat},\times)$ and 
$C^{\bullet}=\mathbb{I}^{\bullet}$. To give a $2$-functor $2_{1}
\rightarrow Coh^{\mathbb{I}^{\bullet}}(\mathcal{A},\mathcal{B})$ is 
to give two $2$-functors $F,G:\mathcal{A}\rightarrow \mathcal{B}$ 
and an object of ${\bf Cat}^{\Delta}(\mathbb{I}^{\bullet},
Y\langle\mathcal{A},\mathcal{B}\rangle(F,G)^{\bullet})$.
By adjunction, the latter data amounts to giving a map 
$\mathbb{I}^{\bullet}\rightarrow Y\langle\mathcal{A},
\mathcal{B}\rangle(F,G)^{\bullet}$
in $\mathbf{Cat}^{\Delta}$. The cosimplicial object 
$Y\langle\mathcal{A},\mathcal{B}\rangle(F,G)^{\bullet}$ 
is described below notation 4.2, thus
$$Y\langle\mathcal{A},\mathcal{B}\rangle(F,G)^{n}=
\begin{cases}
\underset{a\in Ob(\mathcal{A})}\prod \mathcal{B}(Fa,Ga), & \text{if }  n=0\\
\underset{a_{0},...,a_{n}\in Ob(\mathcal{A})}\prod \mathcal{B}(Fa_{0},Ga_{n})
^{\underline{\mathcal{A}}(a_{0},...,a_{n})}, & \text{if } n\geq 1,
\end{cases} $$
where $\underline{\mathcal{A}}(a_{0},...,a_{n})=\mathcal{A}(a_{0},a_{1})
\times...\times \mathcal{A}(a_{n-1},a_{n})$. 
The codegeneracy $s^{0}$ is given by $(u_{a,b})\mapsto (u_{a,a}(1_{a}))$.
Let us compute some cofaces, first 
$$d^{0},d^{1}:Y\langle\mathcal{A},\mathcal{B}\rangle(F,G)^{0}\rightarrow 
Y\langle\mathcal{A},\mathcal{B}\rangle(F,G)^{1}$$
We have, on objects,
$$d^{0}((\alpha_{a}))=(f\mapsto \alpha_{b}\circ Ff) \ {\rm and} 
\ d^{1}((\alpha_{a}))=(f\mapsto Gf\circ \alpha_{a})$$
Next, let's compute $$d^{0},d^{1},d^{2}:Y\langle\mathcal{A},\mathcal{B}\rangle(F,G)^{1}
\rightarrow Y\langle\mathcal{A},\mathcal{B}\rangle(F,G)^{2}$$ We have, on objects,
$$d^{0}((u_{a,b}))=((f:a_{0}\rightarrow a_{1},g:a_{1}\rightarrow a_{2})
\mapsto u_{a_{1},a_{2}}(g)\circ Ff),$$
$$d^{1}((u_{a,b}))=((f:a_{0}\rightarrow a_{1},g:a_{1}\rightarrow a_{2})
\mapsto u_{a_{0},a_{2}}(gf)),$$
$$d^{2}((u_{a,b}))=((f:a_{0}\rightarrow a_{1},g:a_{1}\rightarrow a_{2})
\mapsto Gg\circ u_{a_{0},a_{1}}(f))$$
We appeal now to lemma 5.12 to see what a map $\mathbb{I}^{\bullet}
\rightarrow Y\langle\mathcal{A},\mathcal{B}\rangle(F,G)^{\bullet}$
is. To give a commutative diagram
\[
\xymatrix{
1 \ar[rrr] \ar @<-4pt> [d]_{d^{1}} \ar @<10pt> [d]^{d^{0}} & & &
\underset{a\in Ob(\mathcal{A})}\prod \mathcal{B}(Fa,Ga)
\ar @<-4pt> [d]_{d^{1}} \ar @<10pt> [d]^{d^{0}}\\
\mathbb{I}^{1} \ar @<-1pt> [u]_{p} \ar[rrr]^{f^{1}} & & & 
\underset{a_{0},a_{1}\in Ob(\mathcal{A})}\prod \mathcal{B}(Fa_{0},Ga_{1})
^{\mathcal{A}(a_{0},a_{1})}\ar @<-1pt> [u]_{s^{0}}\\
}
  \]
is to give the data consisting of $(i)$ 1-cells $\alpha_{a}:Fa\rightarrow Ga$ 
for each object $a$ of $\mathcal{A}$, and $(ii)$ a coherence 2-cell 
$\alpha_{f}:Gf\circ \alpha_{a}\rightarrow \alpha_{b}\circ Ff$ 
filling in the square for each $1$-cell $f:a\rightarrow b$, such that 
$\alpha_{1_{a}}:1_{Ga}\circ \alpha_{a}\rightarrow \alpha_{a}\circ 1_{Fa}$ 
is an identity $2$-cell whenever $f=1_{a}$ and for every 
2-cell $\gamma:f\rightarrow g$, $\alpha_{g}(G\gamma \alpha_{a})=
(\alpha_{b} F\gamma)\alpha_{f}$. A map $$f^{2}:\mathbb{I}^{2}\rightarrow 
\underset{a_{0},a_{1},a_{2}\in Ob(\mathcal{A})}\prod \mathcal{B}(Fa_{0},Ga_{2})
^{\underline{\mathcal{A}}(a_{0},...,a_{2})}$$
such that $f^{2}=(d^{0}f^{1},d^{2}f^{1})$ is given on objects by 
$$f^{2}(0)=((f:a_{0}\rightarrow a_{1},g:a_{1}\rightarrow a_{2})
\mapsto G(gf)\circ \alpha_{a_{0}}),$$
$$f^{2}(1)=((f:a_{0}\rightarrow a_{1},g:a_{1}\rightarrow a_{2})
\mapsto Gg\circ \alpha_{a_{1}}\circ Ff),$$
$$f^{2}(2)=((f:a_{0}\rightarrow a_{1},g:a_{1}\rightarrow a_{2})
\mapsto \alpha_{a_{2}}\circ F(gf))$$
To say that $f^{2}c=f^{1}d^{1}$ on arrows is to say
that $(\alpha_{g}Ff)(Gg\alpha_{f})=\alpha_{gf}$.
In conclusion, the $1$-cells of the $2$-category 
$Coh^{\mathbb{I}^{\bullet}}(\mathcal{A},\mathcal{B})$ 
are the quasi-natural (also called lax natural) transformations 
between $2$-functors. A similar argument involving now the 
$2$-functor $2_{[1]}\rightarrow Coh^{\mathbb{I}^{\bullet}}
(\mathcal{A},\mathcal{B})$ shows that
the $2$-cells are the modifications. Therefore 
$\mathcal{A}\boxtimes_{\mathbb{I}^{\bullet}}\mathcal{B}$ 
is Gray's tensor product of $2$-categories. The same considerations 
apply to $Coh^{\mathbb{J}^{\bullet}}(\mathcal{A},\mathcal{B})$.
\\

{\bf Acknowledgements.} This work would not have been possible 
without the help of Andr\'{e} Joyal and Michael Makkai. 
I heartily thank the referees for their comments 
and suggestions and Michael Warren for useful 
discussions related to the material of this article.

\end{document}